# LARGE DEVIATIONS ASSOCIATED WITH POISSON–DIRICHLET DISTRIBUTION AND EWENS SAMPLING FORMULA

BY SHUI FENG [1]

*McMaster University*

Several results of large deviations are obtained for distributions that are associated with the Poisson–Dirichlet distribution and the Ewens sampling formula when the parameter $\theta$ approaches infinity. The motivation for these results comes from a desire of understanding the exact meaning of $\theta$ going to infinity. In terms of the law of large numbers and the central limit theorem, the limiting procedure of $\theta$ going to infinity in a Poisson–Dirichlet distribution corresponds to a finite allele model where the mutation rate per individual is fixed and the number of alleles going to infinity. We call this the finite allele approximation. The first main result of this article is concerned with the relation between this finite allele approximation and the Poisson–Dirichlet distribution in terms of large deviations. Large $\theta$ can also be viewed as a limiting procedure of the effective population size going to infinity. In the second result a comparison is done between the sample size and the effective population size based on the *Ewens sampling formula*.

**1. Introduction.** It is an effective tool to study the infinite dimensional model using their finite dimensional counterpart even though differences, sometimes essential, exist between the two. In this article we focus on a probability distribution, the Poisson–Dirichlet distribution with parameter $\theta > 0$ [henceforth denoted by $PD(\theta)$], on an infinite dimensional space, and its finite dimensional counterpart, the Dirichlet distribution.

Let

$$\nabla = \left\{(p_1, p_2, \ldots): p_1 \geq p_2 \geq \cdots \geq 0, \sum_{k=1}^{\infty} p_k = 1\right\}.$$

Received March 2006; revised March 2007.

[1]Supported in part by the Natural Science and Engineering Research Council of Canada.

*AMS 2000 subject classifications.* Primary 60F10; secondary 92D10.

*Key words and phrases.* Ewens sampling formula, Dirichlet distribution, Poisson–Dirichlet distribution, GEM representation, large deviation.

The infinitely-many-neutral-alleles model is a $\nabla$-valued process describing the evolution of a population under random sampling and parent independent mutation. If the total mutation rate is $\theta$, then the stationary distribution of this process is the $PD(\theta)$. Kingman [15] introduced the $PD(\theta)$ as an asymptotic distribution of the descending order statistics of a symmetric Dirichlet distribution with parameters $K,\alpha$ when $K\to\infty$ and $\alpha\to 0$ in a way such that $\lim_{K\to\infty} K\alpha=\theta$. We use $\mathbf{P}(\theta)=(P_1(\theta),P_2(\theta),\ldots)$ to denote the $\nabla$-valued random variable with distribution $PD(\theta)$. It is equal in law to the normalized jump sizes of Gamma process over the interval $(0,\theta)$ ranked in descending order. A more friendly description of $PD(\theta)$ is as follows.

Consider an i.i.d. sequence of random variables $U_k, k=1,2,\ldots,$ with common distribution $Beta(1,\theta)$. Set

$$X_1=U_1, \qquad X_n=(1-U_1)\cdots(1-U_{n-1})U_n, \qquad n\ge 2. \tag{1.1}$$

Then the law of the descending order statistics of $X_1,X_2,\ldots$ is $PD(\theta)$. The representation in (1.1) is called the GEM representation after R. C. Griffiths, S. Engen and J. W. McCloskey. The sequence $X_k, k=1,2,\ldots,$ corresponds to the size-biased permutation of $PD(\theta)$ and can be obtained through the size-biased sampling of a symmetric Dirichlet distribution with parameters $K,\alpha$ following the procedure of $K\to\infty,\alpha\to 0$, and $K\alpha\to\theta$.

For any fixed $n\ge 1$, let

$$\mathcal{A}_n=\left\{(a_1,\ldots,a_n): a_k\ge 0, k=1,\ldots,n; \sum_{i=1}^n i a_i=n\right\}.$$

Consider a random sample of size $n$ from a Poisson–Dirichlet population and for $k=1,\ldots,n$, define

(1.2) $A_k=$ the number of alleles appearing in the sample exactly $k$ times.

Then $\mathbf{A}_n=(A_1,\ldots,A_n)$ is an $\mathcal{A}_n$-valued random variable with distribution given by the well-known *Ewens sampling formula* [7]:

$$P\{\mathbf{A}_n=(a_1,\ldots,a_n)\}=\frac{n!}{\theta_{(n)}}\prod_{j=1}^n\left(\frac{\theta}{j}\right)^{a_j}\frac{1}{a_j!}, \tag{1.3}$$

where $\theta_{(n)}=\theta(\theta+1)\cdots(\theta+n-1)$.

Consider instead a random sample of size $n$ from a symmetric $Dirichlet(\alpha,\ldots,\alpha)$ distribution with $\theta=K\alpha$, and let $k=\sum_{i=1}^n a_i$, then

$$\begin{aligned} &P\{\mathbf{A}_n=(a_1,\ldots,a_n)\}\\ &\quad=\frac{n!}{\theta_{(n)}}\frac{\alpha^k\Gamma(K+1)}{\Gamma(K-k+1)}\prod_{j=1}^n\left(\frac{\Gamma(j+\alpha)}{\Gamma(j+1)\Gamma(\alpha+1)}\right)^{a_j}\frac{1}{a_j!}. \end{aligned} \tag{1.4}$$

When $K$ goes to infinity, we end up with (1.3).

Hence, many properties of $PD(\theta)$ can be derived from the corresponding results of the Dirichlet distribution with finite alleles. The approximation procedure used is to let $K \to \infty, \alpha \to 0$, and $K\alpha \to \theta$. If $\alpha$ is held constant, then $\theta$ will converge to infinity. This limiting procedure was first suggested by Watterson [22], and some results of the law of large numbers and fluctuations were obtained in Griffiths [11]. But this is not exactly the same as letting $\theta$ go to infinity in $PD(\theta)$ since the latter is a two step limiting procedure: first $K\alpha \to \theta$ as $K \to \infty$, then $\theta$ goes to infinity.

In a large neutral population, the role of mutation is to bring in new type of alleles and to reduce the proportion of existing alleles. The parameter $\theta = 4Nu$ is the population mutation rate with $u$ and $N$ being the individual mutation rate and the effective population size, respectively. The limiting procedure of large $\theta$ is equivalent to a situation where the mutation rate per individual is fixed and the effective population size is large. Motivated by the work of Gillespie [10] on the role of population size in molecular evolution, there have been renewed interests in the asymptotic behaviors of $PD(\theta)$ for large $\theta$ (see [3, 12, 13, 14]).

The first topic in this article is the comparison in terms of large deviations between the finite allele $Dirichlet(\frac{\theta}{K}, \ldots, \frac{\theta}{K})$ distribution and the infinite allele $PD(\theta)$. Two types of limits are considered: the first one is $\theta \to \infty$ followed by $K \to \infty$; the second is $\theta = K \to \infty$. For the first type limit, the following diagram is commutative in terms of the law of large numbers:

$$\begin{array}{ccc} Dirichlet\left(\frac{\theta}{K}, \ldots, \frac{\theta}{K}\right) & \xrightarrow{K \text{ large}} & PD(\theta) \\ \Big\downarrow {\scriptstyle \theta \text{ large}} & & \Big\downarrow {\scriptstyle \theta \text{ large}} \\ \delta_{(1/K, \ldots, 1/K)} & \xrightarrow{K \text{ large}} & \delta_{(0,0,\ldots)}. \end{array} \tag{1.5}$$

In terms of the law of large numbers, the diagram is still commutative when $Dirichlet(\frac{\theta}{K}, \ldots, \frac{\theta}{K})$ and $PD(\theta)$ are replaced by the law of the size-biased sampling and the law of GEM, respectively. But we will show that it fails to be commutative in terms of large deviations. The second type limit corresponds to the diagonal limit in the diagram. The results in [11] and [13] show that the central limit behavior of the finite allele Dirichlet distribution under the second type limit is the same as the infinite allele $PD(\theta)$ for large $\theta$. We will show that the same relation holds for LDP.

The second topic is the comparison between the population and a random sample based on the *Ewens sampling formula*. Consider $\theta$ as a certain population size. Then it would be interesting to see how large the sample size should be so that it behaves the same as the population for large $\theta$.

Here is an outline of the development of this article. In Section 2 the large deviation principles (LDP) are established for the GEM representation and for the finite allele approximations. A detailed comparison is done between the LDPs for $PD(\theta)$ and the GEM, and the corresponding LDPs for Dirichlet distribution and its descending order statistics. An entropy connection is established for the LDP result of $PD(\theta)$ obtained in [3].

LDPs for the conditional and unconditional *Ewens sampling formula* are established in Section 3. The difference between the speed of the two LDPs indicates the strong effect of averaging and finite sample size.

In Section 4 we consider the variable sample size and establish the relation between the sample size $n$, the parameter $\theta$ and the LDP speeds for the total number of alleles, and for the age-class sizes. For the age-class sizes, the sample LDP will be the same as the population LDP if the sample size grows faster than $\theta$.

The reference [4] includes all the terminologies and standard techniques on large deviations used in this article. Since the state spaces encountered here are all compact, we do not make the distinction between a rate function and a good rate function. Generalizations to the two-parameter Poisson–Dirichlet distribution (cf. [17, 18, 19]) will be addressed in a separate article.

**2. LDP associated with $PD(\theta)$.** Let

$$\bar{\nabla} = \left\{(p_1, p_2, \ldots) : p_1 \geq p_2 \geq \cdots \geq 0, \sum_{k=1}^{\infty} p_k \leq 1\right\}$$

be the closure of $\nabla$ equipped with the subspace topology of $R^\infty$.

Let $\mathbf{P} = (P_1, P_2, \ldots)$ be distributed as $PD(\theta)$. In [3] it was shown that an LDP holds for $PD(\theta)$ with speed $\theta$ and rate function

$$I(\mathbf{p}) = \begin{cases} \log \dfrac{1}{1 - \sum_{k=1}^{\infty} p_k}, & (p_1, p_2, \ldots) \in \bar{\nabla}, \displaystyle\sum_{k=1}^{\infty} p_k < 1, \\ \infty, & \text{else.} \end{cases} \tag{2.6}$$

As the first result of this section, we establish the LDP for the GEM. Let

$$\Delta = \left\{(x_1, x_2, \ldots) : x_k \geq 0, k = 1, 2, \ldots; \sum_{i=1}^{\infty} x_i \leq 1\right\},$$

and $\mathbf{X} = (X_1, X_2, \ldots)$ be the GEM in (1.1) generated by the i.i.d. sequence of $Beta(1, \theta)$ random variables $U_1, U_2, \ldots$ . Let $\Pi_\theta^{\text{gem}}$ denote the law of $\mathbf{X}$ on $\Delta$. Since $\mathbf{P}$ and $\mathbf{X}$ stay in different spaces and the ordering operation is not continuous, the LDP for GEM does not follow easily from the LDP for $PD(\theta)$.

For each $n \geq 1$, set

$$\mathbf{X}_n = (X_1, \ldots, X_n),$$

$$\Delta_n = \left\{(x_1, \ldots, x_n) : x_k \geq 0, k = 1, 2, \ldots, n; \sum_{i=1}^{n} x_i \leq 1\right\}$$

and

$$S_n(\mathbf{x}) = \begin{cases} \log \dfrac{1}{1-\sum_{k=1}^n x_k}, & \mathbf{x} \in \Delta_n; \displaystyle\sum_{k=1}^n x_k < 1, \\ \infty, & \text{else.} \end{cases}$$

LEMMA 2.1. *For any $n \geq 1$, let $\Pi_{n,\theta}^{\text{gem}}$ be the law of $\mathbf{X}_n$. Then the family $\{\Pi_{n,\theta}^{\text{gem}} : \theta > 0\}$ satisfies a LDP on $\Delta_n$ with speed $\theta$ and rate function $S_n(\cdot)$.*

PROOF. The LDP follows easily from Lemma 3.1 in [3], the independency of $U_1, \dots, U_n$ and the contraction principle. The rate function has the form

$$\inf\left\{\sum_{i=1}^n \log \frac{1}{1-u_i} : x_k = (1-u_1)\cdots(1-u_{k-1})u_k, k = 1, \dots, n\right\},$$

which is $S_n(\cdot)$ by direct calculation. □

THEOREM 2.2. *The family $\{\Pi_\theta^{\text{gem}} : \theta > 0\}$ satisfies a LDP with speed $\theta$ and rate function*

$$S(\mathbf{x}) = \begin{cases} \log \dfrac{1}{1-\sum_{k=1}^\infty x_k}, & \mathbf{x} = (x_1, x_2, \dots) \in \Delta, \displaystyle\sum_{k=1}^\infty x_k < 1, \\ \infty, & \textit{else.} \end{cases} \tag{2.7}$$

PROOF. First note that the topology on $\Delta$ can be generated by the following metric:

$$d(\mathbf{x}, \mathbf{y}) = \sum_{k=1}^\infty \frac{|x_k - y_k|}{2^k},$$

where $\mathbf{x} = (x_1, x_2, \dots), \mathbf{y} = (y_1, y_2, \dots) \in \Delta$, and the space $\Delta$ is compact. Thus, to establish the LDP, it suffices [20] to verify that

$$\begin{aligned} &\lim_{\delta \to 0} \limsup_{\theta \to \infty} \theta^{-1} \log \Pi_\theta^{\text{gem}}\{d(\mathbf{y}, \mathbf{x}) \leq \delta\} \\ &\quad = \lim_{\delta \to 0} \liminf_{\theta \to \infty} \theta^{-1} \log \Pi_\theta^{\text{gem}}\{d(\mathbf{y}, \mathbf{x}) < \delta\} = -S(\mathbf{x}). \end{aligned} \tag{2.8}$$

For each fixed $\delta_1 > 0, n \geq 1$, and small enough $\delta > 0$, one has

$$\{d(\mathbf{y}, \mathbf{x}) \leq \delta\} \subset \left\{\sup_{1 \leq k \leq n} \{|y_k - x_k|\} \leq \delta_1\right\},$$

which implies that

$$\lim_{\delta\to 0}\limsup_{\theta\to\infty}\theta^{-1}\log\Pi_\theta^{\mathrm{gem}}\{d(\mathbf{y},\mathbf{x})\le\delta\}$$

$$\le \limsup_{\theta\to\infty}\theta^{-1}\log\Pi_{n,\theta}^{\mathrm{gem}}\left\{\sup_{1\le k\le n}\{|y_k-x_k|\}\le\delta_1\right\}$$

$$\le -\inf_{\sup_{1\le k\le n}\{|y_k-x_k|\}\le\delta_1} S_n(y_1,\dots,y_n).$$

Since the very left-hand is independent of $\delta_1$ and $n$, letting $\delta_1$ go to zero and then $n$ go to infinity, we get

$$\lim_{\delta\to 0}\limsup_{\theta\to\infty}\theta^{-1}\log\Pi_\theta^{\mathrm{gem}}\{d(\mathbf{y},\mathbf{x})\le\delta\}\le -S(\mathbf{x}). \tag{2.9}$$

On the other hand, for each fixed $\delta>0$, we can choose $n$ large enough so that

$$\left\{\sup_{1\le k\le n}\{|y_k-x_k|\}<\frac{\delta}{2}\right\}\subset\{d(\mathbf{y},\mathbf{x})<\delta\},$$

which leads to

$$\liminf_{\theta\to\infty}\theta^{-1}\log\Pi_\theta^{\mathrm{gem}}\{d(\mathbf{y},\mathbf{x})<\delta\}$$

$$\ge \liminf_{\theta\to\infty}\theta^{-1}\log\Pi_{n,\theta}^{\mathrm{gem}}\left\{\sup_{1\le k\le n}\{|y_k-x_k|\}<\delta/2\right\}$$

$$\ge -S_n(x_1,\dots,x_n)$$

$$\ge -S(\mathbf{x}).$$

By letting $\delta$ approach zero, it follows that

$$\lim_{\delta\to 0}\liminf_{\theta\to\infty}\theta^{-1}\log\Pi_\theta^{\mathrm{gem}}\{d(\mathbf{y},\mathbf{x})<\delta\}\ge -S(\mathbf{x}), \tag{2.10}$$

which, combined with (2.9), implies the result. □

REMARK. The rate function for GEM has the same form as the rate function for $PD(\theta)$, which is expected because of its exchangeable form.

For each fixed $K\ge 2$, let $(P_1^K,\dots,P_K^K)$ be the decreasing order statistics of a *Dirichlet*$(\frac{\theta}{K},\dots,\frac{\theta}{K})$ random vector. Then the law of $(P_1^K,\dots,P_K^K)$ converges to $PD(\theta)$ in the sense that, for every fixed $r\ge 1$, $(P_1^K,\dots,P_r^K)$ converges to $(P_1,\dots,P_r)$ when $K$ approaches infinity. The LDP for $(P_1^K,\dots,P_K^K)$ when $\theta$ approaches infinity has been established in Theorem 2.1 of [2] with the rate function given by the relative entropy

$$I^K(p_1,\dots,p_K)=\sum_{i=1}^K\frac{1}{K}\log\frac{1/K}{p_i}.$$

We now investigate the structure of diagram (1.5) in terms of the large deviation rate functions.

*Type* I *limit.* $\theta$ goes to infinity followed by $K$ approaches infinity.

THEOREM 2.3.

$$S_r(p_1,\dots,p_r) = \lim_{K\to\infty} \inf\left\{ I^K(p_1,\dots,p_r,q_{r+1},\dots,q_K): p_1\ge\cdots\ge p_r\ge q_{r+1}\ge q_K\ge 0, \sum_{i=1}^r p_i + \sum_{i=r+1}^K q_i = 1\right\}.$$

PROOF. The equality holds trivially if $p_1+\cdots+p_r=1$. We now assume that $\sum_{i=1}^r p_i<1$. It follows that

$$I^K(p_1,\dots,p_r,q_{r+1},\dots,q_K) = \sum_{i=1}^r \frac{1}{K}\log\frac{1}{Kp_i} + \frac{K-r}{K}\log\frac{1}{K} + \frac{1}{K}\log\frac{1}{q_{r+1}\cdots q_K}.$$

Since $\sum_{i=r+1}^K q_i = 1-\sum_{i=1}^r p_i$ and $q_{r+1}\cdots q_K$ reaches the maximum when they are all equal, it follows that

$$\begin{aligned}\inf\{I^K(p_1,\dots,p_r,q_{r+1},\dots,q_K): p_1\ge\cdots\ge p_r\ge q_{r+1}\ge q_K\ge 0,\\ p_1+\cdots+p_r+q_{r+1}+\cdots+q_K=1\}\\ = \frac{r}{K}\log\frac{1}{K} + \frac{1}{K}\sum_{i=1}^r\log\frac{1}{p_i} + \frac{K-r}{K}\log\frac{K-r}{K}\\ + \frac{K-r}{K}\log\frac{1}{1-\sum_{i=1}^r p_i}.\end{aligned}$$

By letting $K$ go to infinity, the equality follows. □

REMARK. The LDP rate function for $(P_1,\dots,P_r)$ under $PD(\theta)$ has been shown to be $S_r(\cdot)$ in [3]. Note that $I^K(\mathbf{q})$ is the relative entropy of $(\frac{1}{K},\dots,\frac{1}{K})$ with respect to $\mathbf{q}$. Thus, we are able to establish certain connections between relative entropy and the LDP for $PD(\theta)$. This also makes a connection between the LDP for *Dirichlet*$(\nu)$ in Dawson and Feng [1, 2] and the LDP for $PD(\theta)$. The LDP speeds for $(P_1,\dots,P_r)$ and $(P_1^K,\dots,P_r^K)$ are both $\theta$.

Let $V_1,\dots,V_{K-1}$ be independent random variables with $V_i$ having a $Beta(\frac{\theta}{K}+1,(K-i)\frac{\theta}{K})$ distribution, and

$$Y_1^K=V_1,\qquad Y_i^K=(1-V_1)\cdots(1-V_{i-1})V_i,\qquad i=2,\dots,K-1, \tag{2.11}$$

be the size-biased sampling of the symmetric $Dirichlet(\frac{\theta}{K},\dots,\frac{\theta}{K})$ distribution.

LEMMA 2.4. *The family of the laws of $(Y_1^K,\dots,Y_{K-1}^K)$ satisfies a LDP on $\Delta_{K-1}$ with speed $\theta$ and rate function*

$$\begin{aligned}S^K(y_1,\dots,y_{K-1})=&\sum_{i=1}^{K-1}\left[\frac{1}{K}\log\frac{1-\sum_{l=1}^{i-1}y_l}{y_i}+\frac{K-i}{K}\log\frac{1-\sum_{l=1}^{i-1}y_l}{1-\sum_{l=1}^{i}y_l}\right]\\&+\sum_{i=1}^{K-1}\left[\frac{1}{K}\log\frac{1}{K+1-i}+\frac{K-i}{K}\log\frac{K-i}{K+1-i}\right].\end{aligned}$$

PROOF. The LDP with speed $\theta$ for $V_i$ for each $i=1,\dots,K-1$ can be established by a direct application of the Laplace method and the rate function for $V_i$ is

$$\begin{aligned}I_{K,i}(v)=&\frac{1}{K}\log\frac{1}{v}+\frac{K-i}{K}\log\frac{1}{1-v}\\&+\frac{1}{K}\log\frac{1}{K+1-i}+\frac{K-i}{K}\log\frac{K-i}{K+1-i},\end{aligned}$$

which implies the result by a combination of independency of $V_1,\dots,V_{K-1}$ and the contraction principle. □

THEOREM 2.5. *The LDP rate function for $(Y_1^K,\dots,Y_r^K)$ for each fixed $r$ does not converge to the LDP rate function for $(X_1,\dots,X_r)$ as $K$ approaches infinity.*

PROOF. It follows from Lemma 2.1 that the LDP rate function for $(X_1,\dots,X_r)$ is $S_r(x_1,\dots,x_r)$, which is finite as long as $\sum_{i=1}^r x_i<1$. By contraction principle, the LDP rate function for $(Y_1^K,\dots,Y_r^K)$ is given by

$$\inf\{S^K(x_1,\dots,x_r,y_{r+1},\dots,y_{K-1}):(x_1,\dots,x_r,y_{r+1},\dots,y_{K-1})\in\Delta_{K-1}\},$$

which is infinite when $x_1=0$ for every $K$. □

Thus, under the Type I limit, the diagram (1.5) is commutative in terms of the LDP rate functions for the ordered distributions but not for the GEM.

*Type* II *limit.* $\theta = K \to \infty$.

When $\theta = K$, the *Dirichlet*$(\frac{\theta}{K}, \dots, \frac{\theta}{K})$ distribution becomes the uniform distribution on the simplex $\{(x_1, \dots, x_K) : x_i \geq 0, i = 1, \dots, K; \sum_{i=1}^{K} x_i = 1\}$. For each fixed $r \geq 1$, the density of $(P_1^K, \dots, P_r^K)$ is given (cf. [21]) by

$$g(p_1, p_2, \dots, p_r) = K(K-1)\cdots(K-r+1)\Gamma(K)L(r, K; B), \tag{2.12}$$

where

$$L(r, K; B) = \int \underset{B}{\dots} \int dx_{r+1} \cdots dx_{K-1},$$

$$B = \Biggl\{(x_{r+1}, \dots, x_{K-1}) : 0 \leq x_i \leq p_r, i = r+1, \dots, K-1; \\ \sum_{i=r+1}^{K-1} x_i \in [1 - a_r - p_r, 1 - a_r]\Biggr\}$$

and $a_r = \sum_{j=1}^{r} p_j$. If $1 - a_r \leq p_r$, then

$$B = \Biggl\{(x_{r+1}, \dots, x_{K-1}) : 0 \leq x_i \leq 1 - p_r, i = r+1, \dots, K-1; \\ \sum_{i=r+1}^{K-1} x_i \in [0, 1 - a_r]\Biggr\}$$

and

$$\begin{aligned} L(r, K; B) &= \int \underset{[0,1-a_r]^{\otimes K-r-2}}{\dots} \int \chi_{\{\sum_{i=r+1}^{K-1} x_i \in [0, 1-a_r]\}} dx_{r+1} \cdots dx_{K-1} \\ &= \frac{(1-a_r)^{K-r-1}}{\Gamma(K-r)}. \end{aligned} \tag{2.13}$$

To deal with the case of $1 - a_r > p_r$, let $X_{r+1}, \dots, X_{K-1}$ be i.i.d. uniform random variables over $(0, p_r)$. Then

$$\begin{aligned} L(r, K; B) &= p_r^{K-r-1} P\Biggl\{1 - a_r - p_r \leq \sum_{i=1}^{K-r-1} X_{r+i} \leq 1 - a_r\Biggr\} \\ &= p_r^{K-r-1}\Biggl[P\Biggl\{\sum_{i=1}^{K-r-1} X_{r+i} \leq 1 - a_r\Biggr\} \\ &\qquad - P\Biggl\{\sum_{i=1}^{K-r-1} X_{r+i} \leq 1 - a_r - p_r\Biggr\}\Biggr]. \end{aligned} \tag{2.14}$$

Let $Z_{r+1},\dots,Z_{K-1}$ be i.i.d. uniform random variables over $(0,1-a_r)$. Since $1-a_r>p_r$, the conditional distribution of $Z_i$ given $Z_i<p_r$ is the same as the law of $X_i$. Hence, by direct calculation, we get

$$
\begin{aligned}
&p_r^{K-r-1}P\Biggl\{\sum_{i=1}^{K-r-1}X_{r+i}\le 1-a_r\Biggr\}\\
&\quad\le(1-a_r)^{K-r-1}P\Biggl\{\sum_{i=1}^{K-r-1}Z_{r+i}\le 1-a_r\Biggr\}\\
&\quad=\frac{(1-a_r)^{K-r-1}}{\Gamma(K-r)}\\
&\quad=\sum_{l=0}^{m}\binom{K-r-1}{l}(1-a_r-p_r)^l p_r^{K-r-1-l}\\
&\qquad\quad\times P\Biggl\{\sum_{i=1}^{K-r-1}Z_{r+i}\le 1-a_r\Bigg|Z_{r+1}>p_r,\dots,Z_{r+l}>p_r,Z_{r+j}\le p_r,\\
&\qquad\qquad\qquad l<j\le K-r-1\Biggr\}\\
&\quad\le C(K,r,m)p_r^{K-r-1-m}P\Biggl\{\sum_{i=m+1}^{K-r-1}X_{r+i}\le 1-a_r\Biggr\},
\end{aligned}
\tag{2.15}
$$

where

$$m=\inf\{k\ge 1: kp_r>1-a_r\}$$

and

$$C(K,r,m)=\sum_{l=0}^{m}\binom{K-r-1}{l}(1-a_r-p_r)^l p_r^{m-l}. \tag{2.16}$$

This combined with (2.13) implies that

$$
\begin{aligned}
\frac{(1-a_r)^{K+m-r-1}}{C(K+m,r,m)\Gamma(K+m-r)}&\le p_r^{K-r-1}P\Biggl\{\sum_{i=1}^{K-r-1}X_{r+i}\le 1-a_r\Biggr\}\\
&\le\frac{(1-a_r)^{K-r-1}}{\Gamma(K-r)}.
\end{aligned}
\tag{2.17}
$$

Similarly, we can prove that

$$p_r^{K-r-1}P\Biggl\{\sum_{i=1}^{K-r-1}X_{r+i}\le 1-a_r-p_r\Biggr\}\le\frac{(1-a_r-p_r)^{K-r-1}}{\Gamma(K-r)}. \tag{2.18}$$

Taking into account (2.14) we get

$$
\begin{aligned}
\frac{(1-a_r)^{K-r-1}}{\Gamma(K-r)} &\geq L(r,K;B)\\
&\geq \frac{(1-a_r)^{K-r-1}}{\Gamma(K-r)} \frac{(1-a_r)^m\Gamma(K-r)}{C(K+m,r,m)\Gamma(K+m-r)}\\
&\quad \times \biggl[1-C(K+m,r,m)\\
&\qquad \times \left(\frac{1-a_r-p_r}{1-a_r}\right)^{K-r-1} \frac{\Gamma(K+m-r)}{(1-a_r)^m\Gamma(K-r)}\biggr].
\end{aligned}
\tag{2.19}
$$

THEOREM 2.6. *The family of the laws of* $(P_1^K,\dots,P_r^K)$ *for each fixed* $r$ *satisfies a LDP with speed* $\theta$ *and rate function* $S_r(p_1,\dots,p_r)$ *as* $\theta=K$ *approaches infinity.*

PROOF. Let

$$\nabla_r=\left\{(q_1,\dots,q_r):0\leq q_r\leq\cdots\leq q_1,\sum_{i=1}^r q_i\leq 1\right\}.$$

For each $\delta>0$ and $(p_1,\dots,p_r)$ in $\nabla_r$, let $B=B((p_1,\dots,p_r),\delta)$ and $\bar{B}=\bar{B}((p_1,\dots,p_r),\delta)$ denote, respectively, the open and closed balls in $\nabla_r$ centered at $(p_1,\dots,p_r)$ with radius $\delta$. Consider the case of $\sum_{k=1}^r p_k<1, p_r>0$ first. By choosing $\delta$ small, we can have $\sum_{k=1}^r q_i<1, q_r>0$ for all $(q_1,\dots,q_r)$ in $\bar{B}$. It then follows from (2.19) that for each $(q_1,\dots,q_r)$ in $\bar{B}$

$$\lim_{\theta\to\infty}\frac{1}{\theta}\log g(q_1,\dots,q_r)=-\log\frac{1}{1-\sum_{k=1}^r q_k}. \tag{2.20}$$

Choose $(q_1^\delta,\dots,q_r^\delta)$ in $\bar{B}$ such that

$$q_r^\delta=\inf\{q_r:(q_1,\dots,q_r)\in\bar{B}\}$$

and

$$\sum_{k=1}^{r-1}q_k^\delta=\sup\left\{\sum_{k=1}^{r-1}q_k:(q_1,\dots,q_r^\delta)\in\bar{B}\right\}.$$

The existence of such point follows from the continuity of the corresponding functions in the above definition.

Since the density function $g(q_1,\dots,q_r)$ is increasing in $q_r$ for fixed $q_1,\dots,q_{r-1}$, and decreasing in $\sum_{k=1}^{r-1}q_k$ for fixed $q_r$, we get that

$$\lim_{\delta\to 0}\liminf_{\theta\to\infty}\frac{1}{\theta}\log P\{(P_1^K,\dots,P_r^K)\in B\}$$

$$
\begin{aligned}
&= \lim_{\delta\to 0}\liminf_{\theta\to\infty}\frac{1}{\theta}\log\int_B g(q_1,\ldots,q_r)\,d\,q_1\cdots d\,q_r\\
&\ge \lim_{\delta\to 0}\liminf_{\theta\to\infty}\frac{1}{\theta}\log\int_B g(q_1,\ldots,q_{r-1},q_r^\delta)\,dq_1\cdots d\,q_r \qquad (2.21)\\
&\ge \lim_{\delta\to 0}\liminf_{\theta\to\infty}\frac{1}{\theta}\log\int_B g(q_1^\delta,\ldots,q_r^\delta)\,dq_1\cdots d\,q_r\\
&= \lim_{\delta\to 0}\log\left(1-\sum_{k=1}^r q_k^\delta\right) = -S_r(p_1,\ldots,p_r),
\end{aligned}
$$

where (2.20) is used to get the second equality.

On the other hand, from (2.19) we have

$$
g(q_1,\ldots,q_r)\le K(K-1)\cdots(K-r+1)\frac{(1-\sum_{k=1}^r q_k)^{K-r-1}\Gamma(K)}{\Gamma(K-r)}.
$$

Thus,

$$
\begin{aligned}
&\lim_{\delta\to 0}\limsup_{\theta\to\infty}\frac{1}{\theta}\log P\{(P_1^K,\ldots,P_r^K)\in\bar{B}\}\\
&\quad\le \lim_{\theta\to\infty}\frac{1}{\theta}\log K(K-1)\cdots(K-r+1)\frac{\Gamma(K)}{\Gamma(K-r)}\\
&\qquad + \lim_{\delta\to 0}\lim_{\theta\to\infty}\frac{1}{\theta}\log\int_{\bar{B}}(1-q_1-\cdots-q_r)^{K-r-1}\,dq_1\cdots d\,q_r \qquad (2.22)\\
&\quad= -\log\frac{1}{1-a_r}.
\end{aligned}
$$

Since the state space is compact, partial LPD holds. From (2.21) and (2.22), all partial rate functions are the same and equal to $S_r(p_1,\ldots,p_r)$ on the set $\{(p_1,\ldots,p_r)\in\nabla_r \sum_{k=1}^r p_k<1, p_r>0\}$. If $\sum_{k=1}^r p_k<1$ and there exists $k_0\le r$ such that $p_k>0$ for $k\le k_0-1$ and $p_k=0$ for $k\ge k_0$, then, for any partial rate function $I'$, we have $I'(p_1,\ldots,p_r)\le S_r(p_1,\ldots,p_r)$ due to the continuity of $S_r(p_1,\ldots,p_r)$ and the lower semi-continuity of $I'$. Thus, the lower bound still holds. If $\sum_{k=1}^r p_k=1$, the lower bound is trivial.

The upper bound holds true in all cases due to the monotonicity of the function $(1-\sum_{k=1}^r q_k)^{K-r-1}$ in $\sum_{k=1}^r q_k$. Thus, (2.21) and (2.22) hold in all cases. This combined with the compactness of the state space implies the theorem. □

For each fixed $r\ge 1$ and the size-biased permutation defined in (2.11), we have the following:

THEOREM 2.7. *The family of the laws of* $(Y_1^K,\ldots,Y_r^K)$ *satisfies a LDP with speed* $\theta$ *and rate function* $S_r(y_1,\ldots,y_r)$ *as* $\theta=K$ *approaches infinity.*

PROOF. Noting that $V_1,\dots,V_r$ are independent and $V_i$ has a $Beta(2,\theta-i)$ distribution, it follows that the law of $(V_1,\dots,V_r)$ satisfies a LDP with speed $\theta$ and rate function

$$\sum_{i=1}^{r}\log\frac{1}{1-v_i}$$

as $\theta=K$ approaches infinity. The theorem then follows easily from the contraction principle. □

REMARK. From Theorems 2.6 and 2.7, we conclude that the LDPs for the finite allele model under the Type II limit are the same for the infinite allele model for large $\theta$. A similar result is expected to hold under the general limit of $\lim_{\theta\to\infty}\frac{\theta}{K}=c>0$.

**3. LDP for Ewens sampling formula.** For each fixed $n\geq 1$, let $\mathbf{A}_n$ be the random partition defined in (1.2). For a given allele proportion $\mathbf{p}=(p_1,p_2,\dots)$ in $\nabla$, and $(a_1,\dots,a_n)$ in $\mathcal{A}_n$, the conditional sampling probability $F_{\mathbf{a}_n}(\mathbf{p})=P\{\mathbf{A}_n=\mathbf{a}_n=(a_1,\dots,a_n)|\mathbf{P}(\theta)=\mathbf{p}\}$ is given by (cf. Kingman [16])

$$F_{\mathbf{a}_n}(\mathbf{p})=C(n,\mathbf{a}_n)\sum p_{l_{11}}p_{l_{12}}\cdots p_{l_{1a_1}}p_{l_{21}}^2p_{l_{22}}^2\cdots p_{l_{2a_2}}^2p_{l_{31}}^3\cdots,\tag{3.23}$$

where

$$C(n,\mathbf{a}_n)=\frac{n!}{\prod_{j=1}^{n}(j!)^{a_j}a_j!}$$

and the summation is over distinct

$$l_{ij},l_{i1}<l_{i2}<\cdots<l_{ia_i},\qquad i=1,\dots,n;j=1,\dots,a_i.$$

If we extend the function $F_{\mathbf{a}_n}(\mathbf{p})$ directly to $\bar{\nabla}$, then we have the following:

LEMMA 3.1. *The function $F_{\mathbf{a}_n}(\mathbf{p})$ is continuous on $\bar{\nabla}$ if and only if $a_1=0$.*

PROOF. Let $k=\sum_{i=1}^{n}a_i$ and $a_1=r$. If $r=k$, then $k=n$. For each $m\geq 1$, let

$$\mathbf{p}_m=\Big(\underbrace{\frac{1}{m},\dots,\frac{1}{m}}_{m},0,\dots\Big)\in\nabla,$$

which converges to $(0,\dots)$ as $m$ goes to infinity. By direct calculation,

$$F_{\mathbf{a}_n}(\mathbf{p}_m)=\binom{m}{n}\left(\frac{1}{m}\right)^n\to\frac{1}{n!}\neq F_{\mathbf{a}_n}((0,\dots)).$$

Next we assume that $1 \le r < k$. For any $m \ge r \vee (k-r)$, let

$$\mathbf{p}_m = \Bigg( \overbrace{\frac{1}{2(k-r)}, \ldots, \frac{1}{2(k-r)}}^{k-r}, \underbrace{\frac{1}{2m}, \ldots, \frac{1}{2m}}_{m}, 0, \ldots \Bigg) \in \nabla,$$

which converges to $\mathbf{q} = (\frac{1}{2(k-r)}, \ldots, \frac{1}{2(k-r)}, 0, \ldots)$ as $m$ goes to infinity. Write

$$F_{\mathbf{a}_n}(\mathbf{p}) = C(n, \mathbf{a}_n) \Bigg[ \sum_1 + \sum_2 \Bigg] p_{l_{11}} p_{l_{12}} \cdots p_{l_{1a_1}} p_{l_{21}}^2 p_{l_{22}}^2 \cdots p_{l_{2a_2}}^2 p_{l_{31}}^3 \cdots,$$

where $\sum_1$ is over indexes such that $\{l_{ij} : i = 2, \ldots, n; j = 1, \ldots, a_i\} = \{1, \ldots, k-r\}$.

By direct calculation, we get

$$\begin{aligned} F_{\mathbf{a}_n}(\mathbf{p}_m) &= C(n, \mathbf{a}_n) \Bigg\{ \binom{k-r}{a_2, \ldots, a_n} \binom{m}{r} \left( \frac{1}{2(k-r)} \right)^{n-r} \left( \frac{1}{2m} \right)^r + o\Big(\frac{1}{m}\Big) \Bigg\} \\ &\to C(n, \mathbf{a}_n) \binom{k-r}{a_2, \ldots, a_n} \left( \frac{1}{2(k-r)} \right)^{n-r} \frac{1}{2^r r!} \neq 0 = F_{\mathbf{a}_n}(\mathbf{q}), \end{aligned}$$

where the $o(\frac{1}{m})$ follows from the fact that each nonzero term in $\sum_2$ involves higher orders of $1/m$.

Thus, $F_{\mathbf{a}_n}(\mathbf{p})$ is not continuous for $a_1 \ge 1$. Next we assume that $a_1 = 0$. Set $N = \max\{l_{2a_2}, \ldots, l_{na_n}\}$. Then for each $M \ge 1$,

$$F_{\mathbf{a}_n}(\mathbf{p}) = F_{\mathbf{a}_n}^1(\mathbf{p}) + F_{\mathbf{a}_n}^2(\mathbf{p}),$$

where

$$\begin{aligned} F_{\mathbf{a}_n}^1(\mathbf{p}) &= C(n, \mathbf{a}_n) \sum_{N \le M} p_{l_{21}}^2 p_{l_{22}}^2 \cdots p_{l_{2a_2}}^2 p_{l_{31}}^3 \cdots, \\ F_{\mathbf{a}_n}^2(\mathbf{p}) &= C(n, \mathbf{a}_n) \sum_{N > M} p_{l_{21}}^2 p_{l_{22}}^2 \cdots p_{l_{2a_2}}^2 p_{l_{31}}^3 \cdots. \end{aligned}$$

Clearly, $F_{\mathbf{a}_n}^1(\mathbf{p})$ is continuous in $\mathbf{p}$. Let $\mathcal{H}$ denote the collection of partitions of $n-1$ obtained from $\mathbf{a}_n$ by removing one individual from the sample.

Then

$$F_{\mathbf{a}_n}^2(\mathbf{p}) \le \frac{C(n, \mathbf{a}_n)}{M} \Bigg\{ \sum_{\mathbf{b} \in \mathcal{H}} \frac{1}{C(n-1, \mathbf{b})} \Bigg\} \to 0 \qquad \text{uniformly in } \mathbf{p}.$$

Thus, for any $\mathbf{p}, \mathbf{q}$ in $\bar{\nabla}$,

$$\begin{aligned} |F_{\{}\mathbf{a}_n(\mathbf{p}) - F_{\mathbf{a}_n}(\mathbf{q})| &\le |F_{\mathbf{a}_n}^1(\mathbf{p}) - F_{\mathbf{a}_n}^1(\mathbf{q})| + |F_{\mathbf{a}_n}^2(\mathbf{p}) - F_{\mathbf{a}_n}^2(\mathbf{q})| \\ &\le |F_{\mathbf{a}_n}^1(\mathbf{p}) - F_{\mathbf{a}_n}^1(\mathbf{q})| + 2\frac{C(n, \mathbf{a}_n)}{M} \Bigg\{ \sum_{\mathbf{b}} \in \mathcal{H}\} \frac{1}{C(n-1, \mathbf{b})} \Bigg\}, \end{aligned}$$

which implies the continuity of $F_{\mathbf{a}_n}(\mathbf{p})$ and the lemma. □

Now Theorem 4.4 in [3] combined with the contraction principle leads to the following:

THEOREM 3.2. *For each* $\mathbf{a}_n$ *in* $\mathcal{A}_n$ *with* $a_1 = 0$, *the family of the laws of* $F_{\mathbf{a}_n}(\mathbf{p})$ *under* $PD(\theta)$ *satisfies an LDP with speed* $\theta$ *and rate function*

$$I_{\mathbf{a}_n}(x) = \inf\left\{\log \frac{1}{1-\sum_{i=1}^{\infty} p_i} : F_{\mathbf{a}_n}(\mathbf{p}) = x\right\}. \tag{3.24}$$

REMARK. If $PD(\theta)$ is replaced by the finite allele symmetric Dirichlet distribution, then the law of $F_{\mathbf{a}_n}(\mathbf{p})$ will satisfy the LDP for all $\mathbf{a}_n$. The function $F_{\mathbf{a}_n}(\mathbf{p})$ can also be extended continuously to $\bar{\nabla}$ by replacing $\sum_{i=1}^{\infty} p_i$ with constant 1. By using this extension, the law of $F_{\mathbf{a}_n}(\mathbf{p})$ satisfies a LDP for all $\mathbf{a}_n$. More detailed discussions on this extension are found in [6].

Next we turn to the large deviations associated with the *Ewens sampling formula*. The state space is now $\mathcal{A}_n$ and the random element is $\mathbf{A}_n$.

THEOREM 3.3. *The family of the laws of* $\mathbf{A}_n$ *under* $PD(\theta)$ *satisfies an LDP with speed* $\log\theta$ *and rate function*

$$I_{\text{esf}}(\mathbf{a}) = n - \sum_{i=1}^{n} a_i. \tag{3.25}$$

PROOF. From the *Ewens sampling formula*, we have

$$\frac{1}{\log\theta}\log P\{\mathbf{A}_n = \mathbf{a}\} = \frac{1}{\log\theta}\left\{\log\frac{n!}{\prod_{j=1}^{n} j^{a_j} a_j!} + \log\frac{\theta^{\sum_{i=1}^{n} a_i}}{\theta_{(n)}}\right\}.$$

The theorem follows by letting $\theta$ go to infinity. □

For the $K$-allele symmetric *Dirichlet*$(\frac{\theta}{K},\dots,\frac{\theta}{K})$ distribution, the sampling formula has the form

$$P\{\mathbf{A}_n = \mathbf{a}\} = \frac{n!}{\prod_{j=1}^{n} a_j!}\frac{K!}{(K-\sum_{i=1}^{n} a_i)!}\frac{1}{\theta_{(n)}}\prod_{j=1}^{n}\left(\frac{\Gamma(j+\theta/K)}{j!\Gamma(\theta/K)}\right)^{a_j}. \tag{3.26}$$

THEOREM 3.4. *Assume that* $\theta$ *grows with* $K$ *and*

$$\lim_{K\to\infty}\frac{\theta}{K} = c \in (0,+\infty].$$

*Then the family of the law of* $\mathbf{A}_n$ *under Dirichlet*$(\frac{\theta}{K},\dots,\frac{\theta}{K})$ *satisfies an LDP with speed* $\log K$ *and rate function* $I_{\text{esf}}(\mathbf{a})$.

Proof. For $c<\infty$,

$$\begin{aligned}\frac{1}{\log K}\log P\{\mathbf{A}_n=\mathbf{a}\} &= \frac{1}{\log K}\Biggl\{\log\frac{n!}{\prod_{j=1}^n j!^{a_j}a_j!}\\ &\qquad +\log\frac{K!}{\theta_{(n)}(K-\sum_{i=1}^n a_i)!}\\ &\qquad +\log\Biggl[\prod_{j=1}^n\biggl(\frac{\Gamma(j+\theta/K)}{\Gamma(\theta/K)}\biggr)^{a_j}\Biggr]\Biggr\}\\ &\to -I_{\mathrm{esf}}(\mathbf{a}) \qquad \text{as } K \text{ goes to infinity.}\end{aligned}$$

For $c=\infty$, by Stirling's formula, we get

$$\begin{aligned}&\frac{1}{\log K}\log P\{\mathbf{A}_n=\mathbf{a}\}\\ &\quad=\frac{1}{\log K}\Biggl\{\log\frac{n!}{\prod_{j=1}^n j!^{a_j}a_j!}\\ &\qquad +\log\frac{K!}{(K-\sum_{i=1}^n a_i)!}\\ &\qquad +\log\Biggl[\frac{1}{\theta_{(n)}}\prod_{j=1}^n\biggl(\frac{\Gamma(j+\theta/K)}{\Gamma(\theta/K)}\biggr)^{a_j}\Biggr]\Biggr\}\\ &\quad=\frac{1}{\log K}\Biggl\{\log\frac{K(K-1)\cdots(K-\sum_{i=1}^n a_i+1)}{K^n}\\ &\qquad +\log\biggl(\frac{\theta^n}{\theta_{(n)}}\biggr)\Biggr\}+O\biggl(\frac{1}{\log K}\biggr)\\ &\quad\to -I_{\mathrm{esf}}(\mathbf{a}) \qquad \text{as } K \text{ goes to infinity.} \qquad \square\end{aligned}$$

Remark. An interesting feature of this theorem is the fact that the large deviation speed and rate function do not depend on the exact speed of $\theta$ as long as it grows no slower than $K$.

**4. Scaling limit with varying sample size.** From the previous section, we see that the LDPs for the conditional *Ewens sampling formula* and the unconditional *Ewens sampling formula* have different speeds. This is due to the averaging and the finite sample size which reduce the randomness. It is thus natural to consider the case of varying sample size.

There are four possible relations between the sample size $n$ and parameter $\theta$, namely:

*Case* A: $n$ fixed and $\theta$ approaches infinity.

*Case* B: $n$ grows with $\theta$ and $\lim_{\theta\to\infty}\theta/n=\infty$.
*Case* C: $n$ grows with $\theta$ and $\lim_{\theta\to\infty}\theta/n=c>0$.
*Case* D: $n$ grows with $\theta$ and $\lim_{\theta\to\infty}\theta/n=0$.
Define

$$\alpha(\theta)=\begin{cases}\log\theta, & \text{Case A},\\ n\log\dfrac{\theta}{n}, & \text{Case B},\\ \theta, & \text{Case C},\\ \theta\log\dfrac{n}{\theta}, & \text{Case D},\end{cases}$$

and

$$\beta(\theta)=\begin{cases}\log\theta, & \text{Case A},\\ \log\dfrac{\theta}{n}, & \text{Case B},\\ \dfrac{\theta}{n}, & \text{Case C},\\ 1, & \text{Case D}.\end{cases}$$

Let

$$K_n(\theta)=\sum_{i=1}^{n}A_i$$

be the total number of distinct alleles in a random sample of size $n$. Let $|S_n^k|$ denote the coefficient of $\theta^k$ in $\theta_{(n)}$. Then the distribution of $K_n(\theta)$ is given by (cf. page 114 in [8])

$$P\{K_n(\theta)=k\}=|S_n^k|\frac{\theta^k}{\theta_{(n)}}. \tag{4.27}$$

The moment generating function of $K_n(\theta)$ is calculated as

$$M(t)=E[e^{tK_n}]=\frac{(e^t\theta)_{(n)}}{\theta_{(n)}}=\frac{\Gamma(\theta e^t+n)\Gamma(\theta)}{\Gamma(\theta e^t)\Gamma(\theta+n)}. \tag{4.28}$$

For large $\theta$, we obtain the following result.

THEOREM 4.1.

$$\Lambda(t)=\lim_{\theta\to\infty}\frac{1}{\alpha(\theta)}\log M(\beta(\theta)t)=\begin{cases}\Lambda_1(t), & \textit{Case A},\\ \Lambda_2(t), & \textit{Case B},\\ \Lambda_3(t), & \textit{Case C},\\ \Lambda_4(t), & \textit{Case D},\end{cases} \tag{4.29}$$

*where,*

$$\Lambda_1(t)=\begin{cases}nt, & \textit{if } t>-1,\\ (t+1)-n, & \textit{else},\end{cases} \tag{4.30}$$

$$\Lambda_2(t) = \begin{cases} t, & \textit{if } t > -1, \\ -1, & \textit{else}, \end{cases} \tag{4.31}$$

$$\Lambda_3(t) = \frac{1}{c}\{[c\log c - (1+c)\log(1+c)] \tag{4.32}$$
$$+ [(1+ce^{ct})\log(1+ce^{ct}) - ce^{ct}\log ce^{ct}]\}$$

$$\Lambda_4(t) = e^t - 1. \tag{4.33}$$

Proof. Case A and Case C follow easily from direct calculations. For Case B, we use the Stirling formula to get

$$\begin{aligned} &\frac{1}{\alpha(\theta)}\log M(\beta(\theta)t) \\ &\quad = \frac{1}{\log\theta/n}\left[\log\left(\frac{1+(\theta/n)^{t+1}}{1+\theta/n}\right) + \log\left(1+\frac{1}{(\theta/n)^{t+1}}\right)^{(\theta/n)^{t+1}} \right. \\ &\qquad\qquad \left. - \log\left(1+\frac{1}{(\theta/n)}\right)^{(\theta/n)}\right] + o(1) \\ &\quad \to \begin{cases} t, & \text{if } t > -1, \\ -1, & \text{else.} \end{cases} \end{aligned}$$

For Case D, using Stirling's formula again, we get

$$\begin{aligned} \frac{1}{\alpha(\theta)}\log M(\beta(\theta)t) &= e^t\frac{\log(1+n/\theta e^{-t})}{\log n/\theta} - \frac{\log(1+n/\theta)}{\log n/\theta} + o(1) \\ &\to e^t - 1. \end{aligned}$$

□

Theorem 4.2. *In Case* A, *the family of the laws of* $K_n(\theta)$ *on space* $\{1,\dots,n\}$ *under* $PD(\theta)$ *satisfies an LDP with speed* $\log\theta$ *and rate function*

$$I(k) = n - k.$$

Proof. For each $k = 1,\dots,n$, it follows from (4.27) that

$$\begin{aligned} \lim_{\theta\to\infty}\frac{1}{\log\theta}\log P\{K_n(\theta)=k\} &= \lim_{\theta\to\infty}\frac{1}{\log\theta}\log\frac{\theta^k}{\theta_{(n)}} \\ &= -(n-k), \end{aligned}$$

which implies the result. □

Theorem 4.3. *In Case* B, *the family of the laws of* $K_n(\theta)/n$ *on space* $[0,1]$ *under* $PD(\theta)$ *satisfies an LDP with speed* $n\log\frac{\theta}{n}$ *and rate function*

$$I(x) = 1 - x.$$

Proof. As the coefficient of $\theta^k$ in $\theta_{(n)}$, the number $|S_n^k|$ lies between $\frac{(n-1)!}{(k-1)!}$ and $\binom{n-1}{k-1}\frac{(n-1)!}{(k-1)!}$. It follows from Stirling's formula that

$$\begin{aligned}\log \frac{(n-1)!\theta^k}{(k-1)!\theta_{(n)}} &= n\biggl[-\log\left(1+\frac{\theta}{n}\right)+\frac{k}{n}\left(\log\frac{\theta}{n}-\log\frac{k}{n}\right)\\ &\qquad -\log\left(1+\frac{1}{\theta/n}\right)^{\theta/n}+o(1)\biggr]\\ &= -\alpha(\theta)\left[1-\frac{k}{n}+o(1)\right]\end{aligned}$$

and

$$\begin{aligned}\log\binom{n-1}{k-1} &= -(n-1)\left[\frac{k-1}{n-1}\log\frac{k-1}{n-1}+\left(1-\frac{k-1}{n-1}\right)\log\frac{n-k}{n-1}+o(1)\right]\\ &= -\alpha(\theta)o(1).\end{aligned}$$

Thus, it follows from (4.27) that

$$\frac{1}{\alpha(\theta)}\log P\left\{\frac{K_n(\theta)}{n}=\frac{k}{n}\right\}=-\left(1-\frac{k}{n}+o(1)\right).$$

For each $x$ in $[0,1]$ and $\delta>0$, the total number of integers in

$$\{1,\dots,n\}\cap[nx-n\delta,nx+n\delta]$$

is of the magnitude of $n\delta$ which is $\alpha(\theta)o(1)$. Hence,

$$\begin{aligned}&\lim_{\delta\to0}\liminf_{\theta\to\infty}\frac{1}{\alpha(\theta)}\log P\left\{\left|\frac{K_n(\theta)}{n}-x\right|<\delta\right\}\\ &\quad=\lim_{\delta\to0}\limsup_{\theta\to\infty}\frac{1}{\alpha(\theta)}\log P\left\{\left|\frac{K_n(\theta)}{n}-x\right|\le\delta\right\}\\ &\quad=-(1-x),\end{aligned}$$

which combined with the compactness of $[0,1]$ implies the result. □

Theorem 4.4. *In Case* C, *the family of the laws of* $K_n(\theta)/n$ *on space* $[0,1]$ *under* $PD(\theta)$ *satisfies an LDP with speed* $\theta$ *and rate function*

$$I(x)=\sup_{t\in R}\{tx-\Lambda_3(t)\}.$$

Proof. Extend the law of $K_n(\theta)/n$ to the whole real line $R$ and denote the extension by $P_{\theta,n}$.

According to Definition 2.3.5 in [4], the function $\Lambda_3(t)$ is essentially smooth. By the Gärtner–Ellis Theorem, a LDP holds for the family $\{P_{\theta,n} : \theta > 0, n = 1, \ldots\}$ with rate function

$$\tilde{I}(x) = \sup_{t \in R}\{tx - \Lambda_3(t)\}.$$

To prove the theorem, it suffices to show that $\tilde{I}(x) = \infty$ for $x \notin [0,1]$. Since

$$\lim_{t \to -\infty} \Lambda_3(t) = \frac{1}{c}[c \log c - (1+c)\log(1+c)],$$

by letting $t$ approach negative infinity, one gets that $\tilde{I}(x) = \infty$ for $x < 0$. By direct calculation,

$$\begin{aligned}
&(1 + ce^{ct})\log(1 + ce^{ct}) - ce^{ct}\log(ce^{ct}) \\
&\quad = \log(1 + ce^{ct}) + (ce^{ct})\log\left(1 + \frac{1}{ce^{ct}}\right) \\
&\quad \approx ct \qquad \text{as } t \text{ approaches positive infinity.}
\end{aligned}$$

Hence, by choosing $t > 0$ large, we get that $\tilde{I}(x) = \infty$ for $x > 1$. □

Since $\Lambda_4(t)$ is also essentially smooth, by an argument similar to that used in the proof of Theorem 4.4, we get the following:

Theorem 4.5. *In Case* D, *the family of the laws of* $K_n(\theta)/\theta \log \frac{n}{\theta}$ *on space* $[0, +\infty)$ *under* $PD(\theta)$ *satisfies an LDP with speed* $\theta \log \frac{n}{\theta}$ *and rate function*

$$I(x) = x \log x - x + 1.$$

As a by-product of these LDPs, we get the following weak laws of large numbers.

Corollary 4.1.

$$\text{Case A}, \quad \lim_{\theta \to \infty} K_n(\theta) = n, \tag{4.34}$$

$$\text{Case B}, \quad \lim_{\theta \to \infty} \frac{K_n(\theta)}{n} = 1, \tag{4.35}$$

$$\text{Case C}, \quad \lim_{\theta \to \infty} \frac{K_n(\theta)}{n} = \log\left(1 + \frac{1}{c}\right)^c, \tag{4.36}$$

$$\text{Case D}, \quad \lim_{\theta \to \infty} \frac{K_n(\theta)}{\theta \log(n/\theta)} = 1. \tag{4.37}$$

PROOF. In Case A, the rate function has a unique zero point at $n$. Thus, (4.34) holds. Similarly, results (4.35) and (4.37) follow from the fact that the corresponding rate functions have unique zero point at 1. To prove (4.37), rewrite $\Lambda_3(t)$ as

$$\Lambda_3(t) = \frac{1}{c}[F(ce^{ct}) - F(c)],$$

where $F(x) = (1+x)\log(1+x) - x\log(x)$. Then, by a change of variable of $u = ce^{ct}$, it follows that

$$I_3(x) = \frac{1}{c}\left[\sup_{u>0}\{x\log u - F(u)\} - (x\log c - F(c))\right]. \tag{4.38}$$

For each fixed $x \in (0,1)$, the supremum in (4.38) is achieved at a unique $u_x$ satisfying

$$e^x = \left(1 + \frac{1}{u}\right)^u.$$

The cases of $x = 0$ or 1 correspond to $u \to 0$ or $u \to \infty$, respectively. Thus, $I_3(x) > 0$ unless $u_x = c$ or, equivalently, for $x = \log(1 + \frac{1}{c})^c$. This leads to (4.36). □

REMARK. It follows from (4.37) that the LDP in Case D is very similar to the case studied in [9] where $\theta$ is fixed and $n$ approaches infinity.

The parameter $\theta$ is proportional to certain effective population size when the individual mutation rate is held constant. The limiting procedure of $\theta$ approaching infinity is thus equivalent to letting the population size go to infinity. This next result compares sample size $n$ with the effective population size through the study of the age-class size. It reveals how large the sample size should be so that the sample and the effective population will behave the same as $\theta$ approaches infinity.

Let $X_{1,n}, \ldots, X_{n,n}$ be the age-class sizes in the sample. Then from Donnelly and Tavaré [5], one has

$$P\{X_{1,n} = k\} = \frac{\theta}{n}\frac{\binom{n}{k}}{\binom{\theta+n-1}{k}} = \frac{\theta}{n}\frac{n!}{(n-k)!}\frac{(\theta+n-k-1)!}{(\theta+n-1)!}. \tag{4.39}$$

THEOREM 4.6. *In Case* A, *the family of the laws of* $X_{1,n}$ *on space* $\{1,\ldots,n\}$ *satisfies a LDP with speed* $\log\theta$ *and rate function* $k-1$.

PROOF. A direct calculation gives

$$\begin{aligned}&\lim_{\theta\to\infty}\frac{1}{\log\theta}\log P\{X_{1,n}=k\}\\ &\quad=1-\sum_{i=1}^{k}\lim_{\theta\to\infty}\frac{\log(\theta+n-i)}{\log\theta}=-(k-1).\end{aligned}\tag{4.40}$$

□

THEOREM 4.7. *The family of the laws of* $\frac{X_{1,n}}{n}$ *on space* $E=[0,1]$ *satisfies a LDP with speed* $\gamma(\theta)$ *and rate function* $I(x)$ *given respectively by*

$$\gamma(\theta)=\begin{cases} n\log\dfrac{\theta}{n}, & \textit{Case B},\\ \theta, & \textit{Case C},\\ \theta, & \textit{Case D},\end{cases}\tag{4.41}$$

*and*

$$I(x)=\begin{cases} x, & \textit{Case B},\\ I_c(x), & \textit{Case C},\\ \log\dfrac{1}{1-x}, & \textit{Case D},\end{cases}\tag{4.42}$$

*where*

$$\begin{aligned}I_c(x)&=\frac{1}{c}[(c+1)\log(c+1)\\ &\qquad+(1-x)\log(1-x)-(c+1-x)\log(c+1-x)].\end{aligned}\tag{4.43}$$

REMARK 1. The last case is the same as the LDP for the GEM. In other words, when $n$ grows faster than $\theta$, the random sample behaves the same as the population.

PROOF OF THEOREM 4.7. It follows from Stirling's formula that

$$\begin{aligned}&\log P\{X_{1,n}=k\}\\ &\quad=\log\left[\frac{(n+1)^{n+1/2}(n+\theta-k)^{\theta+n-k-1/2}}{(n-k+1)^{n-k+1/2}(\theta+n)^{\theta+n-1/2}}\right]\\ &\qquad+\log\frac{\theta}{n}+o(1)\\ &\quad=I_1+I_2+I_3+I_4+I_5+I_6+o(1),\end{aligned}\tag{4.44}$$

where

$$I_1=\log\frac{\theta}{n},$$

$$I_2 = \frac{1}{2}\left[\log\frac{n+1}{n} + \log\left(\frac{\theta}{n}+1\right)\right],$$

$$I_3 = \left(k - \frac{1}{2}\right)\log\left(1 + \frac{1}{n} - \frac{k}{n}\right),$$

$$I_4 = -\left(k + \frac{1}{2}\right)\log\left(1 + \frac{\theta}{n} - \frac{k}{n}\right),$$

$$I_5 = n\log\frac{n+1}{n+1-k},$$

$$I_6 = (\theta+n)\log\left(1 - \frac{k}{\theta+n}\right).$$

For each $x$ in $E$, let $\lfloor nx \rfloor$ denote the integer part of $nx$. It is not hard to see that in Case B $I_1 + I_2 + I_3 + I_5 + I_6 = o(\gamma(\theta))$ uniformly in $k/n$ and $I_4(\frac{\lfloor nx \rfloor}{n})/\gamma(\theta)$ converges to $-x$ uniformly in $x$. In Case C, $I_1 + I_2 = o(\gamma(\theta))$ uniformly in $k/n$ and $\frac{1}{\gamma(\theta)}[I_3 + I_4 + I_5 + I_6](\frac{\lfloor nx \rfloor}{n}) \to -I_c(x)$ uniformly in $x$. Thus, we have

$$\lim_{\theta\to\infty}\frac{1}{\gamma(\theta)}\log P\{X_{1,n} = \lfloor nx \rfloor\} = -I(x). \tag{4.45}$$

For each $x$ in $E$ and $\delta > 0$, choose $n$ large enough so that $\lfloor nx \rfloor$ is in the interval $(nx - n\delta, nx + n\delta)$. Then

$$\begin{aligned} P\{X_{1,n} = \lfloor nx \rfloor\} &\le P\{|X_{1,n}/n - x| < \delta\} \\ &\le P\{|X_{1,n}/n - x| \le \delta\} \\ &= \sum_{|k/n - x| \le \delta} P\{X_{1,n} = k\} \\ &\le 2(n\delta + 1)\max_{k\in[nx-n\delta, nx+n\delta]} P\{X_{1,n} = k\}. \end{aligned} \tag{4.46}$$

By the uniform convergence, and (4.46), we get

$$\begin{aligned} &\lim_{\delta\to 0}\lim_{\theta\to\infty}\frac{1}{\gamma(\theta)}\log P\{|X_{1,n}/n - x| < \delta\} \\ &\quad = \lim_{\delta\to 0}\lim_{\theta\to\infty}\frac{1}{\gamma(\theta)}\log P\{|X_{1,n}/n - x| \le \delta\} = -I(x). \end{aligned} \tag{4.47}$$

The maneuvering in Case D is a little different. By a reorganization of the terms, we have

$$\log P\{X_{1,n} = k\} = J_1 + J_2 + J_3 + J_4 + J_5 + o(1), \tag{4.48}$$

where

$$\begin{aligned}
J_1 &= \log\theta - \log n,\\
J_2 &= \frac{1}{2}\left[\log\left(\frac{1+1/n}{1+1/n-k/n}\right) + \log\left(\frac{1+\theta/n}{1+\theta/n-k/n}\right)\right],\\
J_3 &= k\log\frac{1+1/n-k/n}{1+\theta/n-k/n},\\
J_4 &= n\log\frac{(n+1)(\theta+n-k)}{(n+1-k)(\theta+n)},\\
J_5 &= \theta\log\frac{\theta/n+1-k/n}{\theta/n+1}.
\end{aligned}$$

Noting that

$$\begin{aligned}
&n\delta \min_{k\in[nx-n\delta,nx+n\delta]} P\{X_{1,n}=k\}\\
&\quad\le P\{|X_{1,n}/n - x|\le\delta\}\\
&\quad\le 2(n\delta+1)\max_{k\in[nx-n\delta,nx+n\delta]} P\{X_{1,n}=k\},
\end{aligned}\tag{4.49}$$

by taking the logarithm, the term $-\log n$ in $J_1$ is cancelled by $\log(n\delta)$ and the term $\log\theta$ clearly grows slower than $\gamma(\theta)$.

First consider the case $0\le x<1$, and choose $\delta$ small enough so that $x+\delta<1$.

It is clear that $J_2=o(\theta)$, and $\frac{1}{\theta}J_5$ converges uniformly to $\log(1-y)$ over $[x-\delta,x+\delta]$ as $\theta\to\infty, k/n\to y$. For $J_3$ and $J_4$, we have

$$\begin{aligned}
\frac{1}{\theta}J_3 &= \frac{k}{\theta}\log\left(1+\frac{(1-\theta)/n}{1+\theta/n-k/n}\right)\\
&= \frac{k}{n}\frac{n}{\theta}\log\left(1+\frac{(1-\theta)/n}{1+\theta/n-k/n}\right)\\
&\to -\frac{y}{1-y} \qquad \text{as } \theta\to\infty, k/n\to y
\end{aligned}$$

and

$$\begin{aligned}
\frac{1}{\theta}J_4 &= \frac{n}{\theta}\log\left(1+\frac{(\theta-1)k}{n^2+n\theta+(1-k)(n+\theta)}\right)\\
&= \frac{n}{\theta}\log\left(1+\frac{(\theta-1)k/n^2}{1+\theta/n+(1/n-k/n)(1+\theta/n)}\right)\\
&\to \frac{y}{1-y} \qquad \text{as } \theta\to\infty, k/n\to y.
\end{aligned}$$

Summing-up all the terms, one gets

$$\begin{aligned}\lim_{\delta\to 0}\lim_{\theta\to\infty}\frac{1}{\theta}\log P\{|X_{1,n}/n-x|<\delta\}\\ &=\lim_{\delta\to 0}\lim_{\theta\to\infty}\frac{1}{\theta}\log P\{|X_{1,n}/n-x|\le\delta\}\\ &=-\log\frac{1}{1-x}.\end{aligned}\tag{4.50}$$

Finally, for $x=1$, the result still holds from the fact that $P\{X_{1,n}=k\}$ is decreasing in $k$ and

$$\max_{k\in[nx-n\delta,nx+n\delta]}P\{X_{1,n}=k\}\le P\{X_{1,n}=\lfloor n(1-\delta)\rfloor\}.$$ □

This result can be generalized to the first $r$ family sizes in a sample of size $n$. From Donnelly and Tavaré [5], one has, for $k_i\ge 1, i=1,\dots,r$,

$$\begin{aligned}P\{X_{1,n}=k_1,\dots,X_{r,n}=k_r\}\\ &=\frac{(\theta/n)^r}{(1-k_1/n)\cdots(1-k_1/n-\cdots-k_{r-1}/n)}\\ &\quad\times\frac{n!}{(n-k_1-\cdots-k_r)!}\frac{(\theta+n-k_1-\cdots-k_r-1)!}{(\theta+n-1)!}.\end{aligned}\tag{4.51}$$

This formula is very similar to (4.39) and all asymptotic calculations are almost the same except one needs to replace $k$ there with $k_1+\cdots+k_r$. The state space is

$$\left\{(k_1,\dots,k_r):k_i\ge 1, i=1,\dots,r;\sum_{j=1}^r k_j\le n\right\}$$

in Case A and

$$\left\{(x_1,\dots,x_r):x_i\in[0,1], i=1,\dots,r;\sum_{j=1}^r x_j\le 1\right\}$$

TABLE 1

| **Case** | **Speed** | **Rate function** |
|---|---|---|
| A | $\log\theta$ | $\sum_{i=1}^r k_i-r$ |
| B | $n\log\frac{\theta}{n}$ | $\sum_{i=1}^r x_i$ |
| C | $\theta$ | $I_c(\sum_{i=1}^r x_i)$ |
| D | $\theta$ | $\log\frac{1}{1-\sum_{i=1}^r x_i}$ |

in all other cases. The corresponding LDP results are summarized in Table 1.

**Acknowledgments.** I wish to thank Paul Joyce for informing me of the results in [5], and helpful discussions. I also wish to thank the referee for valuable comments and suggestions.

Department of Mathematics and Statistics
McMaster University
Hamilton, Ontario
Canada L8S 4K1
E-mail: shuifeng@mcmaster.ca